\documentclass[11pt,a4paper]{article}
\usepackage[utf8]{inputenc}
\usepackage{amsmath}
\usepackage{amsfonts}
\usepackage{amssymb}
\usepackage{amsthm}
\usepackage{setspace}
\usepackage{hyperref}
\usepackage{graphicx}
\usepackage{stmaryrd}
\usepackage{mathtools}
\usepackage{subcaption}
\usepackage[T1]{fontenc}
\usepackage{xcolor}
\usepackage{float}
\usepackage{tikz-cd}
\usepackage[noadjust]{cite}
\usepackage{cite}
\usepackage{ragged2e}
\usepackage[margin=1.2in]{geometry}
\graphicspath{ {./images/} }
\newtheorem{thm}{Theorem}[section]
\newtheorem{prop}[thm]{Proposition}
\newtheorem{lemma}[thm]{Lemma}
\newtheorem{coro}[thm]{Corollary}
\newtheorem*{thmA}{Theorem A}
\newtheorem*{thmB}{Theorem B}

\theoremstyle{definition}
\newtheorem{defi}[thm]{Definition}

\newtheorem{rem}[thm]{Remark}

\newcommand{\quotient}[2]{{\raisebox{.2em}{$#1$}\left/\raisebox{-.2em}{$#2$}\right.}}

\begin{document}

\begin{center}
{\large \textbf{Automorphisms of large-type free-of-infinity Artin groups}}
\end{center}

\begin{center}
Nicolas Vaskou
\end{center}

\begin{abstract}
\centering \justifying We compute explicitly the automorphism and outer automorphism group of all large-type free-of-infinity Artin groups. Our strategy involves reconstructing the associated Deligne complexes in a purely algebraic manner, i.e. in a way that is independent from the choice of standard generators for the groups.
\end{abstract}

\noindent \rule{7em}{.4pt}\par

\small



\normalsize

\section{Introduction}

Artin groups form a large family of groups that have drawn an increasing attention in the past few decades. They are defined as follows. Start from a simplicial graph $\Gamma$ with finite vertex set $V(\Gamma)$ and finite edge set $E(\Gamma)$. For every edge $e^{ab}\in E(\Gamma)$ with vertices $a, b \in V(\Gamma)$, associate an integer coefficient (or label) $m_{ab} \geq 2$. Then $\Gamma$ is the \textbf{presentation graph} of an \textbf{Artin group} $A_{\Gamma}$ whose presentation is the following:
$$A_{\Gamma} \coloneqq \langle \ V(\Gamma) \ | \ \underbrace{aba \cdots}_{m_{ab}} = \underbrace{bab \cdots}_{m_{ab}} \text{ for every } e^{ab} \in E(\Gamma) \ \rangle.$$
We will say that $m_{ab} = \infty$ whenever $a$ and $b$ are not adjacent. The elements of $V(\Gamma)$ are called the \textbf{standard generators} of $A_{\Gamma}$, and their number $|V(\Gamma)|$ is called the \textbf{rank} of $A_{\Gamma}$. We suppose throughout the paper that $\Gamma$ is connected, a weak condition that ensures the group does not trivially decompose as a free product of infinite groups.

Artin groups are cousins of Coxeter groups: whenever $\Gamma$ defines an Artin group $A_{\Gamma}$, it also defines a \textbf{Coxeter group} $W_{\Gamma}$ that can be obtained from the presentation of $A_{\Gamma}$ by adding the relation $a^2 = 1$ for every standard generator. While Coxeter groups are rather well-understood, much less is known about Artin groups in general. Although they are conjectured to have many properties with various flavours (torsion-free-ness, solvable conjugacy problem, $K(\pi, 1)$-conjecture, biautomaticity, acylindrical hyperbolicity, CAT($0$)-ness, etc.), proving any property in full generality has remained exceedingly complicated. In general, such conjectures are solved for more specific classes of Artin groups, assuming additional properties about their presentation graphs (\cite{charney1995k}, \cite{vaskou2021acylindrical}, \cite{huang2019metric}, \cite{haettel2019xxl}).

In \cite{goldsborough2023random}, the authors introduced a notion of randomness for presentation graphs. They showed that for each of the aforementioned conjectures, there is a class for which it has been solved that has a “non-trivial” asymptotic size within the family of all Artin groups. However, there are two (intrinsically related) questions regarding Artin groups that have remained more mysterious: that of solving the isomorphism problem, and that of computing their automorphism groups. The \textbf{isomorphism problem} for Artin groups asks what can be said about two presentations graphs $\Gamma$ and $\Gamma'$ assuming their corresponding Artin groups $A_{\Gamma}$ and $A_{\Gamma'}$ are isomorphic. In \cite{vaskou2023isomorphism}, the author solved this problem for \textbf{large-type} Artin groups (those with coefficients at least $3$). The study of isomorphisms between Artin groups is inherent to the study of automorphisms of Artin groups. As for the isomorphism problem, the study of the automorphisms of Artin groups has turned out to be quite difficult. The most famous results are that of right-angled Artin groups (\cite{droms1987isomorphisms}, \cite{servatius1989automorphisms}, \cite{laurence1992automorphisms}). The situation becomes even more complicated when we introduce non-commuting relations. The only results on Artin groups that are not right-angled concern the class of “connected large-type triangle-free” Artin groups introduced by Crisp (\cite{crisp2005automorphisms}, \cite{an2022automorphism}).

The main goal of the present paper is to compute the automorphism group and outer automorphism group of a larger family of Artin groups. This provides a first example of a family that does not have trivial asymptotic size for which this problem has been solved.

Before giving our main results, we introduce some terminology. An Artin group $A_{\Gamma}$ is said to be \textbf{free-of-infinity} if $m_{ab} < \infty$ for every pair of standard generators $a, b \in V(\Gamma)$. The class of free-of-infinity Artin groups is particularly interesting: in \cite{godelle2012basic}, Godelle and Paris proved that many important conjectures of Artin groups can be solved in full generality if they were solved for free-of-infinity Artin groups.
\medskip

Our main result is the following:

\begin{thmA} Let $A_{\Gamma}$ be a large-type free-of-infinity Artin group of rank at least $3$. Then $Aut(A_{\Gamma})$ is generated by the following automorphisms:
\medskip

\noindent \textbf{1. Conjugations.} Also known as inner automorphisms, these are the automorphisms of the form $\varphi_g : h \mapsto g h g^{-1}$ for some $g \in A_{\Gamma}$.
\medskip

\noindent \textbf{2. Graph automorphisms.} Every label-preserving graph automorphism $\phi \in Aut(\Gamma)$ induces a permutation of the standard generators, hence an automorphism of the group.
\medskip

\noindent \textbf{3. The global inversion.} This is the order $2$ automorphism $\iota$ that sends every standard generator to its inverse.
\medskip

In particular, $Out(A_{\Gamma})$ is finite and isomorphic to $Aut(\Gamma) \times \left(\quotient{\mathbf{Z}}{2 \mathbf{Z}}\right)$.
\end{thmA}

Note that it would not be possible to extend the previous theorem to all large-type Artin groups, as these contain a fourth type of automorphisms called “Dehn twists automorphisms” (see \cite{crisp2005automorphisms}).

Before explaining our strategy we recall a few notions related to Artin groups. If $A_{\Gamma}$ is an Artin group, every induced subgraph $\Gamma' \subseteq \Gamma$ generates a subgroup $\langle V(\Gamma') \rangle \subseteq A_{\Gamma}$. It is a well-known result that this subgroup is isomorphic to the Artin group $A_{\Gamma'}$ itself (\cite{van1983homotopy}). Such subgroups are called \textbf{standard parabolic subgroups}, and their conjugates are called \textbf{parabolic subgroups}. If $\Gamma'$ is such that the associated Coxeter group $W_{\Gamma'}$ is finite, then $A_{\Gamma'}$ is called \textbf{spherical}. In large-type Artin groups, a parabolic subgroup $g A_{\Gamma'} g^{-1}$ is spherical if and only if $\Gamma'$ is a single vertex, a single edge, or empty.

In \cite{charney1995k}, Charney and Davis introduced a combinatorial complex $X_{\Gamma}$ known as the \textbf{Deligne complex}, which has since then proved multiple times to be an incredibly efficient geometric tool in the study of Artin groups. This complex is constructed from the combinatorics of the spherical parabolic subgroups (see Definition \ref{DefiDeligne}). For large-type Artin groups, its geometry is generally better understood.

We now come back to our strategy for proving Theorem A. The Deligne complex $X_{\Gamma}$ is a priori very much dependent on the choice of presentation graph $\Gamma$ for the associated Artin group. That said, if we find a way to reconstruct $X_{\Gamma}$ with purely algebraic objects, then any automorphism of the Artin group will preserve the structure of these algebraic objects, and hence preserve the Deligne complex itself. This approach allows to build an action of the automorphism group $Aut(A_{\Gamma})$ on the Deligne complex, from which we can recover a full description of $Aut(A_{\Gamma})$. This kind of technique was originally used by Ivanov (\cite{ivanov2002mapping}) to study the automorphisms of mapping class groups and has since then been extended to other groups like Higman's group (\cite{martin2017cubical}) or graph products of groups (\cite{genevois2018automorphism}).

In our case, when the Artin groups considered are large-type and free-of-infinity, we find a way to “reconstruct” the associated Deligne complexes in a purely algebraic manner. We obtain the following:

\begin{thmB} Let $A_{\Gamma}$ be a large-type free-of-infinity Artin group of rank at least $3$. Then $X_{\Gamma}$ can be reconstructed in a way that is invariant under isomorphisms (in particular, independent of $\Gamma$). Consequently, there is a natural combinatorial action of $Aut(A_{\Gamma})$ on $X_{\Gamma}$, and this action can be described explicitly.
\end{thmB} 

For the precise description of the action, see Theorem \ref{ThmIsomorphicD}. We briefly mention that the above theorem directly provides a solution to the isomorphism problem for large-type free-of-infinity Artin groups (see Theorem \ref{TheoremRigidity}). We do not include this result as one of our main theorems as it is already a consequence of \cite[Theorem B]{vaskou2023isomorphism}.

Once again, Theorem B cannot be extended to all large-type Artin groups. Indeed, in this class, two Artin groups may be isomorphic while their presentation graphs are not (see \cite{brady2002rigidity}). In particular, one can prove that their associated Deligne complexes are not isomorphic either.

Consider a large-type Artin group $A_{\Gamma}$. A first step to reconstruct the associated Deligne complex $X_{\Gamma}$ is to reconstruct what are called the “type $2$” vertices of the complex. These vertices are in one-to-one correspondence with the spherical parabolic subgroups of $A_{\Gamma}$ on $2$ generators, which have been proved in \cite{vaskou2023isomorphism} to be invariant under automorphisms.

A second step, that is the main technical result of this paper, is to be able to reconstruct the “type $1$” vertices of the Deligne complex. This is where the hypothesis of being free-of-infinity comes into play. The stabilisers of these type $1$ vertices are parabolic subgroups of type $1$ of $A_{\Gamma}$, which have been proved to not be preserved under automorphisms for general large-type Artin groups (\cite[Theorem H]{vaskou2023isomorphism}). However, they are preserved for large-type free-of-infinity Artin groups. That said, this correspondence between the type $1$ vertices and the parabolic subgroups of type $1$ of $A_{\Gamma}$ is far from being a bijection, as infinitely many type $1$ vertices may lie on a common “standard tree” and hence have the same stabiliser. On such a standard tree, it can generally be quite hard to give a purely algebraic condition that translates when two type $1$ vertices should be “adjacent”. Such a condition actually cannot exist for large-type Artin groups in general. In this paper we provide such a condition, under the hypothesis that the groups are also free-of-infinity.

At this point, we will have reconstructed (most of) the $1$-skeleton of the Deligne complex in a purely algebraic way. With a little bit more work, we will then be able to reconstruct the whole complex, proving Theorem B.

\bigskip
\textbf{Organisation of the paper:} Section 2 serves as a preliminary section where we introduce various algebraic and geometric notions, such as the definition of the Deligne complex. In Section 3, we focus on large-type free-of-infinity Artin groups, and we reconstruct their Deligne complexes purely algebraically, proving Theorem B. Finally in Section 4, we use this algebraic description of the Deligne complex to prove Theorem A.

\section{The Deligne complex}

In this section we introduce the Deligne complex and various related geometric objects. We start by defining this complex:

\begin{defi} \label{DefiDeligne}
Let $A_{\Gamma}$ be any Artin group. The \textbf{Deligne complex} $X_{\Gamma}$ associated with $A_{\Gamma}$ is the simplicial complex constructed as follows:
\begin{itemize}
\itemsep0em 
\item The vertices of $X_{\Gamma}$ are the left-cosets $g A_{\Gamma'}$, where $g \in A_{\Gamma}$ and $A_{\Gamma'}$ is any spherical standard parabolic subgroup of $A_{\Gamma}$.
\item Every string of inclusion of the form $g_0 A_{\Gamma_0} \subsetneq \cdots \subsetneq g_n A_{\Gamma_n}$
spans an $n$-simplex.
\end{itemize}
The group $A_{\Gamma}$ acts on $X_{\Gamma}$ by left multiplication.
\end{defi}

The geometry of the Deligne complex for large-type (and more generally $2$-dimensional) Artin groups is better understood than the general case. This is mostly due to the fact that it can be given a CAT(0) metric:

\begin{defi} \label{DefiMetric}
Let $A_{\Gamma}$ be a large-type Artin group. The \textbf{Moussong metric} on $X_{\Gamma}$ is the piecewise Euclidean metric $d$ that is defined by requiring every $2$-simplex of the form $\{1\} \subsetneq \langle a \rangle \subsetneq A_{ab}$ to be isometric to the unique Euclidean triangle satisfying the following (up to isometry):
$$d(\{1\}, \langle a \rangle) \coloneqq 1, \ \ \ \ \angle_{\langle a \rangle}(\{1\}, A_{ab}) \coloneqq \frac{\pi}{2}, \ \ \ \ \angle_{A_{ab}}(\{1\}, \langle a \rangle) \coloneqq \frac{\pi}{2 m_{ab}}.$$
\end{defi}

\begin{thm} \label{ThmXCAT(0)} \textbf{(\cite[Proposition 4.4.5]{charney1995k})} Let $A_{\Gamma}$ be a large-type Artin group. Then the Deligne complex $X_{\Gamma}$ is $2$-dimensional and the Moussong metric is \textup{CAT(0)}.
\end{thm}

\begin{rem} \label{RemGammaBar}
The fundamental domain of the action of $A_{\Gamma}$ on $X_{\Gamma}$ is the subcomplex $K_{\Gamma}$ whose vertices are the spherical standard parabolic subgroups of $A_{\Gamma}$. Since $\{1\}$ is contained in every spherical parabolic subgroup, the corresponding vertex is attached to every spherical standard parabolic subgroup $A_{\Gamma'}$. In particular, $K_{\Gamma}$ is a cone whose apex is $\{1\}$. It is not hard to see that the boundary of $K_{\Gamma}$ is graph-isomorphic to the barycentric subdivision $\Gamma_{bar}$ of $\Gamma$. Hence we will often write $\Gamma_{bar}$ to denote the boundary of $K_{\Gamma}$.
\end{rem}

\begin{figure}[H]
\centering
\includegraphics[scale=0.78]{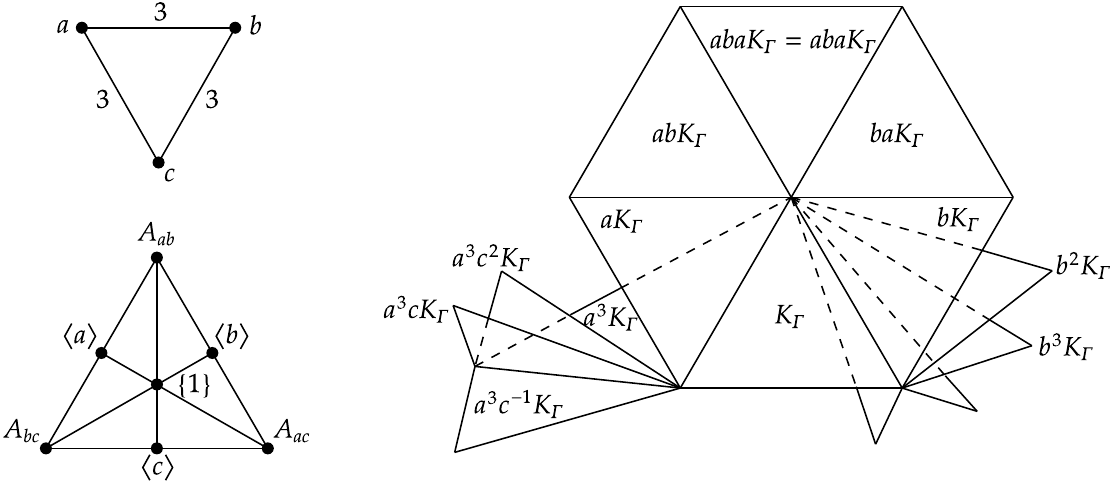}
\caption{\underline{Top-left:} A graph $\Gamma$ defining a large-type Artin group. \underline{Bottom-left:} The fundamental domain $K_{\Gamma}$ of $X_{\Gamma}$. \underline{Right:} the Deligne complex $X_{\Gamma}$. For drawing purposes, only the edges contained in $X_{\Gamma}^{(1)-ess}$ have been drawn.}
\label{FigDeligneComplex}
\end{figure}

In Section 3, the first step to reconstruct the Deligne complex purely algebraically will be to reconstruct the subcomplex corresponding to the points with non-trivial stabiliser. This subcomplex, that we define thereafter, plays an important role throughout the paper:

\begin{defi}
The \textbf{essential 1-skeleton} of $X_{\Gamma}$ is the subcomplex of the $1$-skeleton $X_{\Gamma}^{(1)}$ defined by
$$X_{\Gamma}^{(1)-ess} \coloneqq \bigcup\limits_{g \in A_{\Gamma}} g \Gamma_{bar}.$$
\end{defi}

\begin{rem} \label{RemX1essCone}
Since $K_{\Gamma}$ is the cone-off of $\Gamma_{bar}$, the Deligne complex $X_{\Gamma}$ can be obtained from $X_{\Gamma}^{(1)-ess}$ by coning-off the translates $g \Gamma_{bar}$, for all $g \in A_{\Gamma}$.
\end{rem}

At last, we want to talk about fixed-point sets and types of elements:

\begin{defi} The \textbf{fixed set} of a subset $S \subseteq A_{\Gamma}$ for the action on $X_{\Gamma}$ is the set
$$Fix(S) \coloneqq \{p \in X_{\Gamma} \ | \ \forall g \in S, \ g \cdot p = p \}.$$
\end{defi}

\begin{defi} \label{DefiType} The \textbf{type} of a parabolic subgroup $g A_{\Gamma'} g^{-1}$ is the integer $|V(\Gamma')|$. The \textbf{type} of a simplex of $\sigma \subseteq X_{\Gamma}$ is the type of its stabiliser $G_{\sigma}$.
\end{defi}

\begin{rem}
If $P = g A_{\Gamma'} g^{-1} = h A_{\Gamma''} h^{-1}$, then it follows from \cite[Theorem A]{vaskou2023isomorphism} that $|V(\Gamma')| = |V(\Gamma'')|$. In particular, the type of a parabolic subgroup is always well-defined.
\end{rem}

\begin{lemma} \label{LemmaClassificationByType} \textbf{(\cite[Lemma 8]{crisp2005automorphisms})}
Let $A_{\Gamma}$ be a large-type Artin group, and let $A_{\Gamma'} \subseteq A_{\Gamma}$ be a standard parabolic subgroup. Then:
\\ $\bullet$ $type(A_{\Gamma'}) \geq 3$ or $V(\Gamma') = \{a, b\}$ for some $a, b \in V(\Gamma)$ with $m_{ab} = \infty$ $\Longleftrightarrow$ $Fix(A_{\Gamma'}) = \emptyset$. 
\\ $\bullet$ $\Gamma' = e^{ab}$ for some $a, b \in V(\Gamma)$ with $m_{ab} < \infty$ $\Longleftrightarrow$ $Fix(A_{\Gamma'})$ is the type $2$ vertex $A_{ab}$.
\\ $\bullet$ $\Gamma' = \{a\}$ $\Longleftrightarrow$ $Fix(A_{\Gamma'})$ is a tree called the \textbf{standard tree} of $a$.
\\ The same applies to all parabolic subgroups, as $Fix(g A_{\Gamma'} g^{-1}) = g Fix(A_{\Gamma'})$.
\end{lemma}

\begin{rem}
When $\Gamma$ is connected, the standard tree $Fix(\langle a \rangle)$ contains infinitely many type $1$ vertices, including the vertex $\langle a \rangle$ itself. It also contains infinitely many type $2$ vertices (unless $a$ lies at the tip of an even-labelled leaf - see (\cite[Lemma 2.17]{vaskou2023isomorphism})).
\end{rem}

\section{Reconstructing the Deligne complex algebraically}

Let $A_{\Gamma}$ be a large-type free-of-infinity Artin group. That is, every pair of distinct standard generators $a, b \in V(\Gamma)$ has a coefficient $3 \leq m_{ab} < \infty$. Note that the free-of-infinity condition forces $\Gamma$ to be a complete graph.

This section is dedicated to reconstructing the Deligne complex of $A_{\Gamma}$ in a purely algebraic way. This will allow to build a suitable action of $Aut(A_{\Gamma})$ onto $X_{\Gamma}$, proving Theorem B.
\medskip

\begin{figure}[ht]
\centering
\includegraphics[scale=0.9]{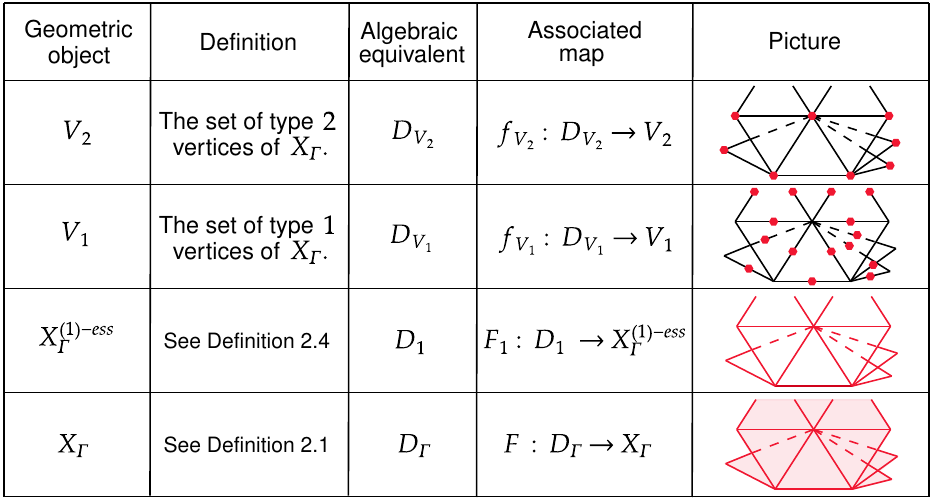}
\caption{Notations used throughout Section 3.}
\label{FigureNotations}
\end{figure}

\noindent \textbf{Strategy and notation:} Our strategy can be divided in four steps. At each step, the goal will be to introduce a set of algebraic objects that “corresponds” to a set of geometric objects of $X_{\Gamma}$. These various correspondences will be made explicit through maps that will be bijections, graph isomorphisms or combinatorial isomorphisms, depending on the context. We sum up the various notations that will be used in Figure \ref{FigureNotations}.
\medskip

Recall that $A_{\Gamma}$ is any large-type free-of-infinity Artin group of rank at least $3$. Our first goal will be to construct algebraic equivalent of the sets $V_2$ and $V_1$ of vertices of type $2$ and type $1$ respectively. Then, we will describe when the algebraic objects corresponding to the elements of $V_2$ and $V_1$ should be “adjacent”, allowing to reconstruct the type $1$ edges (i.e. $X_{\Gamma}^{(1)-ess}$). We start with the following definition:

\begin{defi} \label{DefiDV2} Let $D_{V_2}$ be the set of type $2$ spherical parabolic subgroups of $A_{\Gamma}$.
\end{defi}

\noindent Note that while parabolic subgroups are in general not purely algebraically defined, in (\cite[Corollary 4.17]{vaskou2023isomorphism}) the author proved that type $2$ spherical parabolic subgroups of large-type Artin groups are in fact invariant under automorphisms. In particular, the set $D_{V_2}$ is also invariant under automorphisms.

\begin{lemma} \label{LemmaType2VertAlg} The map $f_{V_2} : D_{V_2} \rightarrow V_2$ defined as follows is a bijection:
\\(1) For every subgroup $H \in D_{V_2}$, $f_{V_2}(H)$ is the fixed set $Fix(H)$;
\\(2) For every vertex $v \in V_2$, $f_{V_2}^{-1}(v)$ is the local group $G_v$.
\end{lemma}

\noindent \textbf{Proof:} This directly follows from Lemma \ref{LemmaClassificationByType}.
\hfill\(\Box\)
\bigskip

\noindent Reconstructing the type $1$ vertices of $X_{\Gamma}$ algebraically will be harder, because showing that parabolic subgroups of type $1$ are invariant under automorphisms does not imply a result similar to Lemma \ref{LemmaType2VertAlg} for the type $1$ vertices of $X_{\Gamma}$. We start by introducing the following property:

\begin{defi} \label{DefiEdgeProp} A pair of subgroups $(H_1, H_2) \in D_{V_2} \times D_{V_2}$ is said to have the \textbf{adjacency property} if there exists a subgroup $H_3 \in D_{V_2}$ such that we have
\begin{align*}
&(A1) \ H_i \cap H_j \neq \{1\}, \ \forall i, j \in \{1,2,3\}; \\
&(A2) \ \bigcap\limits_{i=1}^3 H_i = \{1\}.
\end{align*}
\end{defi}

\noindent Definition \ref{DefiEdgeProp} really is geometric in essence, as highlighted in the next lemma.

\begin{lemma} \label{LemmaEdgePropGeom} A pair $(H_1, H_2)$ has the \text{adjacency property} relatively to a third subgroup $H_3$ if and only if the following hold:
\\(1) The three $H_i$'s are distinct subgroups.
\\(2) The three intersections $(H_i \cap H_j)$ are parabolic subgroups of type $1$, and they are distinct. Equivalently, the sets $Fix(H_i \cap H_j)$ are distinct standard trees.
\\(3) The standard trees $Fix(H_i \cap H_j)$ intersect each other pairwise, but the triple-intersection is trivial.
\end{lemma}

\noindent \textbf{Proof:} $(\Rightarrow)$ Suppose that $(H_1, H_2)$ has the adjacency property relatively to a third subgroup $H_3$. Let $i, j, k \in \{1,2,3\}$ be distinct, and suppose that $H_i = H_j$. Then
$$\{1\} \overset{(A2)} = H_i \cap H_j \cap H_k = H_i \cap H_k \overset{(A1)} \neq \{1\},$$
a contradiction. This proves $(1)$.

In particular, any intersection $H_i \cap H_j$ is a proper non-trivial intersection of parabolic subgroups of type $2$ of $A_{\Gamma}$, hence is a parabolic subgroup of type $1$ of $A_{\Gamma}$, by (\cite[Proposition 2.22.(1)]{vaskou2023isomorphism}). It follows that each $Fix(H_i \cap H_j)$ is a standard tree. This proves $(2)$.

Finally, on one hand the three standard trees intersect each other pairwise, as for instance the intersection of $Fix(H_i\cap H_j)$ and $Fix(H_i \cap H_k)$ is the vertex $Fix(H_i)$. On the other hand, the intersection of the three standard trees is the intersection of all the pairwise intersections. It is trivial because the three vertices $Fix(H_i)$, $Fix(H_j)$ and $Fix(H_k)$ are distinct, as their corresponding subgroups are. This proves $(3)$. 

$(\Leftarrow)$ Suppose that the three subgroups $H_1, H_2, H_3 \in D_{V_2}$ satisfy the properties $(1)$, $(2)$ and $(3)$ of the lemma. The fact that all the intersections $(H_i \cap H_j)$ are parabolic subgroups of type $1$ directly implies $(A1)$.

The subgroups $H_i \cap H_j$ and $H_i \cap H_k$ are parabolic subgroups of type $1$ of $A_{\Gamma}$, so their intersection is a parabolic subgroup of $A_{\Gamma}$ as well, by (\cite[Proposition 2.22.(1)]{vaskou2023isomorphism}). By (\cite[Proposition 2.22.(4)]{vaskou2023isomorphism}), this intersection cannot be a parabolic subgroup of type $1$ of $A_{\Gamma}$, because $H_i \cap H_j$ and $H_i \cap H_k$ are distinct. So it must be trivial. This imples $(A2)$.
\hfill\(\Box\)

\begin{prop} \label{PropAdjEqAdj} Consider two subgroups $H_1, H_2 \in D_{V_2}$. The following are equivalent:
\\(1) The two type $2$ vertices $v_1, v_2$ of $X_{\Gamma}$ defined by $v_i \coloneqq f_{V_2}(H_i)$ are at combinatorial distance $2$ in $X_{\Gamma}^{(1)-ess}$.
\\(2) The pair $(H_1, H_2)$ satisfies the adjacency property.
\end{prop}

\noindent Note that the minimal combinatorial distance one can have between two type $2$ vertices of $X_{\Gamma}^{(1)-ess}$ is $2$, so the previous proposition gives an algebraic description of when two type $2$ vertices of $X_{\Gamma}$ are “as close as possible”.
\medskip

Our main tool to prove the proposition will be the technical Lemma \ref{LemmaTrianglesInFundamentalDomains}. To prove it, we will need the following definitions, leading up to Theorem \ref{ThmGB}. Most can be found in \cite{mccammond2002fans}, although we sometimes state them in a slightly more specific case which serves our narrative.

\begin{defi}
A \textbf{disc diagram} $D$ is a non-empty finite contractible $2$-dimensional polygonal complex that embeds in a disc. $D$ is called \textbf{non-singular} if it is homeomorphic to a disc. If $D$ is non-singular, then it has a \textbf{boundary} that we denote by $\partial D$, and an \textbf{interior} that we denote by $int(D)$, which respectively are the image of the boundary and the interior of the unit disc $\mathcal{D}^2$ under the previous homeomorphism $\mathcal{D}^2 \rightarrow D$.
\end{defi}

\begin{defi}
Let $D$ be a non-singular disc diagram, and suppose that the restriction of the attaching map $\mathcal{D}^2 \rightarrow D$ to any polygon is injective. A \textbf{corner} of $D$ is a pair $(v, f)$ where $f$ is a polygon of $D$ and $v$ is a vertex of $f$.

We will denote by $Corners(v)$ the set of corners whose first component is $v$, and by $Corners(f)$ the set of corners whose second component is $f$.

A non-singular disc diagram is said to be \textbf{angled} if every corner $c$ is given a real number $\angle c$ called the \textbf{angle} at $c$.
\end{defi}

\begin{defi}
Let $D$ be a non-singular angled disc diagram. We will denote by $D_0$ the set of vertices of $D$, and by $D_2$ the set of polygons of $D$. The \textbf{curvature} of an element of $D_0$ or of $D_2$ is defined as follows:
\begin{align*}
&\forall v \in int(D_0), \ curv(v) \coloneqq 2 \pi - \left( \sum\limits_{c \in Corners(v)} \angle c \right), \\
&\forall v \in \partial D_0, \ curv(v) \coloneqq \pi - \left( \sum\limits_{c \in Corners(v)} \angle c \right), \\
&\forall f \in D_2, \ curv(f) \coloneqq 2 \pi - \left( \sum\limits_{c \in Corners(f)} (\pi - \angle c) \right).
\end{align*}
\end{defi}

\begin{thm} \label{ThmGB} \textbf{(\cite[Theorem 4.6]{mccammond2002fans}, Combinatorial Gauss-Bonnet)} Let $D$ be a non-singular angled disc diagram. Then we have
$$\sum\limits_{v \in D_0} curv(v) + \sum\limits_{f \in D_2} curv(f) = 2 \pi.$$
\end{thm}

\begin{lemma} \textbf{(Triangle of standard trees)} \label{LemmaTrianglesInFundamentalDomains}
Let $A_{\Gamma}$ be a large-type Artin group. Let $v_1$, $v_2$ and $v_3$ be three distinct type $2$ vertices of $X_{\Gamma}$, and suppose that the three geodesics connecting the vertices, with respect to the \textup{CAT(0)} metric, are contained in distinct standard trees that intersect pairwise but whose triple intersection is empty. Then the triangle formed by these three geodesics is contained in a single translate $g K_{\Gamma}$ of the fundamental domain.
\end{lemma}

\noindent \textbf{Proof:} Let $T$ be the closed path in $X_{\Gamma}$ obtained by concatenating the three geodesics $\gamma_{ij}$ connecting $v_i$ and $v_j$ for $\{i, j\} \subseteq \{1, 2, 3\}$. By hypothesis, any of the $\gamma_{ij}$'s is contained in a standard tree, which ensures that $T$ is contained in $X_{\Gamma}^{(1)-ess}$. In particular, $T$ is a concatenation of edges. Using that $X_{\Gamma}$ is CAT($0$) (Theorem \ref{ThmXCAT(0)}), we can thus apply \cite[Theorem A.1]{bader2024cat} to conclude that $T$ bounds a unique embedded disc $D$ in $X_{\Gamma}$. Note that $D$ is a non-singular disc diagram, whose boundary is $\partial D = T$ and whose interior we'll denote by $int(D)$. We want to prove that $D$ is contained in a single fundamental domain $g K_{\Gamma}$. To do so we suppose that this is not the case, and we will exhibit a contradiction. We want to apply the Gauss-Bonnet formula on $D$. By construction, $D$ is a combinatorial subcomplex of $X_{\Gamma}$ whose simplices are triangles corresponding to inclusions of the form $\{g\} \subsetneq g \langle a \rangle \subsetneq g A_{ab}$. To make the use of the Gauss-Bonnet formula easier, we decide to see $D$ with a coarser combinatorial structure: the one obtained by removing every edge of type $0$ and every vertex of type $0$ in $D$. We equip $D$ with the metric $\left. d \right|_D$ where $d$ is the usual Moussong metric on $X_{\Gamma}$ (see Definition \ref{DefiMetric}). Note that the attaching maps of polygons are injective. With this structure, we see $D$ as an angled disc diagram in the obvious way: for any corner $c = (v, f)$ where $v$ is a type $2$ vertex, the local group at $v$ is a dihedral Artin group with coefficient $m_c \geq 3$, and we let $\angle c \coloneqq \pi / m_c$. If $v$ is type $1$, then $\angle c \coloneqq \pi$. In particular, type $1$ vertices do not contribute to curvature. By Theorem \ref{ThmGB}, we have
$$\sum\limits_{\text{type } 2 \text{ vertices } v \text{ in } D} curv(v) \ \ + \sum\limits_{\text{polygons } f \text{ in } D} curv(f) = 2 \pi. \ \ \ (*)$$
We rewrite this in a manner that is easier to deal with. Let $D_2^i$ be the set of polygons in $D$ that don't contain any element of $\{v_1, v_2, v_3 \}$, $D_2^v$ be the set of polygons in $D$ that contain at least one of $v_1$, $v_2$ or $v_3$, $D_0^i$ be the set of type $2$ vertices in $int(D)$, $D_0^b$ be the set of type $2$ vertices of $\partial D \backslash \{v_1,v_2,v_3\}$, and $D_0^v$ be the set $\{v_1,v_2,v_3\}$. Then:
\\ $\bullet$ Let $C_2^i \coloneqq \sum\limits_{f \in D_2^i} curv(f)$. Consider a polygon $f \in D_2$, and for every $c \in Corners(f)$ let $m_c$ be as before. Then
$$curv(f) = 2 \pi - \left( \sum\limits_{c \in Corners(f)} \left(\pi - \frac{\pi}{m_c} \right) \right).$$
Note that $m_c \geq 3$ for all $c \in Corners(f)$, so eventually $\pi - \frac{\pi}{m_c} \geq \frac{2 \pi}{3}$. Because $Corners(f)$ contains at least $3$ elements, we obtain
$$curv(f) \leq 2 \pi - 3 \cdot \left(\frac{2 \pi}{3} \right) = 0.$$
It follows that $C_2^i \leq 0$ as well. Note that as soon as one polygon has at least $4$ edges, or as soon as the coefficient of one of the local groups is at least $4$, we have $curv(f) < 0$ and thus $C_2^i < 0$.
\\ $\bullet$ Let $C_0^i \coloneqq \sum\limits_{v \in D_0^i} curv(v)$. Because $X_{\Gamma}$ is CAT($0$), the systole of the link of any vertex $v$ in $X_{\Gamma}$ is at least $2 \pi$. In particular, if $v \in D_0^i$, the systole of the link of $v$ in $M$ is at least $2 \pi$. Since $D$ is embedded in $X_{\Gamma}$, it follows that the sum of the angles around $v$ in $D$ is at least $2 \pi$. In particular, $curv(v) \leq 0$ and thus $C_0^i \leq 0$.
\\ $\bullet$ Let $C_0^b \coloneqq \sum\limits_{v \in D_0^b} curv(v)$. Any $v \in D_0^b$ belongs to a side of $T$ that is a geodesic, so its angle with $D$ must satisfy $\angle_v D \geq \pi$. It follows that $curv(v) = \pi - \angle_v D \leq 0$, and thus $C_0^b \leq 0$ as well.
\\ $\bullet$ Let $C_0^v \coloneqq \sum\limits_{v_i \in D_0^v} curv(v_i)$ and let $C_2^v = \sum\limits_{f \in D_2^v} curv(f)$. Any $v_i$ belongs to $\lambda_i \geq 1$ polygons of $D$. By construction of the Deligne complex, the angle $\angle_{v_i} D$ is precisely $\lambda_i \cdot \frac{\pi}{m_i}$, where $m_i \geq 3$ is the coefficient of the local group at $v_i$. Each of the $\lambda_i$ polygons $f$ of $D$ containing $v_i$ is such that
\begin{align*}
curv(f) =\ &2 \pi - (\pi - \angle (v_i, f)) - \left( \sum\limits_{c \in Corners(f) \backslash \{(v_i, f)\} } \left(\pi - \frac{\pi}{m_c} \right) \right) \\
\overset{(**)} \leq &2 \pi - \left(\pi - \frac{\pi}{m_i} \right) - 2 \cdot \left(\pi - \frac{\pi}{3} \right) \\
\leq \ &\frac{\pi}{m_i} - \frac{\pi}{3}.
\end{align*}
The inequality $(**)$ comes from the fact that $Corners(f)$ has at least $2$ other elements than $(v_i, f)$ and that the angle of any $c \in Corners(f)$ is at most $\pi/3$ because every local group has coefficient at least $3$. Note that if $f$ has at least $4$ edges then we obtain a strict inequality $curv(f) < \frac{\pi}{m_i} - \frac{\pi}{3}$. Summing everything, we obtain
\begin{align*}
C_0^v + C_2^v = &\sum\limits_{v_i \in D_0^v} curv(v_i) + \sum\limits_{f \in D_2^v} curv(f) \\
\overset{(***)} \leq &\sum\limits_{ i \in \{1, 2, 3 \} } \left(\pi - \lambda_i \cdot \frac{\pi}{m_i} \right) + \sum\limits_{ i \in \{1, 2, 3 \} } \lambda_i \cdot \left(\frac{\pi}{m_i} - \frac{\pi}{3} \right) \\
= \ &3 \pi - \sum\limits_{ i \in \{1, 2, 3 \} } \lambda_i \cdot \frac{\pi}{3} \leq 2 \pi.
\end{align*}
Note that it is easy to check that the inequality $(***)$ holds no matter if the polygons containing the $v_i$'s are distinct or if there are polygons of $D$ that contain several of the $v_i$'s. We now notice two things. The first is that as soon as one of the $v_i$'s is contained inside two distinct polygons of $D$, then $\lambda_i \geq 2$ and $C_0^v + C_2^v < 2 \pi$. The second is that if a polygon containing one of the $v_i$'s has at least $4$ edges, then $curv(f) < \pi/m_i - \pi/3$ and thus $C_0^v +C_2^v < 2 \pi$ as well.

\medskip
\noindent With this setting, the equation $(*)$ becomes:
$$C_2^i+C_0^i+C_0^b+(C_0^v+C_2^v) = 2\pi.$$
Note that this equation can hold only if the four terms on the left-hand side are maximal, i.e.:
\\ $\bullet$ $C_2^i = 0$. In particular, every polygon $f$ in $D_2^i$ is a triangle, whose corners $(v, f)$ are such that the local group at $v$ has coefficient exactly $3$.
\\ $\bullet$ $C_0^i =0$. In particular, the sum of the angles around any vertex of $D_0^i$ is exactly $2 \pi$.
\\ $\bullet$ $C_0^b = 0$, i.e. the angles along the sides of $T$ are exactly $\pi$.
\\ $\bullet$ $C_0^v+C_2^v = 2 \pi$. In particular, each of the $v_i$'s is contained in a single polygon of $D$, which is always a triangle.
\medskip

By hypothesis $D$ does not contain a single polygon, and it is not hard to see that in that case there must be polygons in $D$ that do not contain any of the $v_i$'s (in other words, $D_2^i$ is non-trivial). The first of the above four points implies that every polygon in $D_2^i$ is an Euclidean equilateral triangle. We then use the following:
\medskip

\noindent \underline{Claim:} If $f_1$ and $f_2$ are two triangles of $D$ that share an edge, then $f_2$ is a translate of $f_1$. In particular, any two triangles in $D$ are translates of each other's.
\medskip

\begin{figure}[H]
\centering
\includegraphics[scale=0.7]{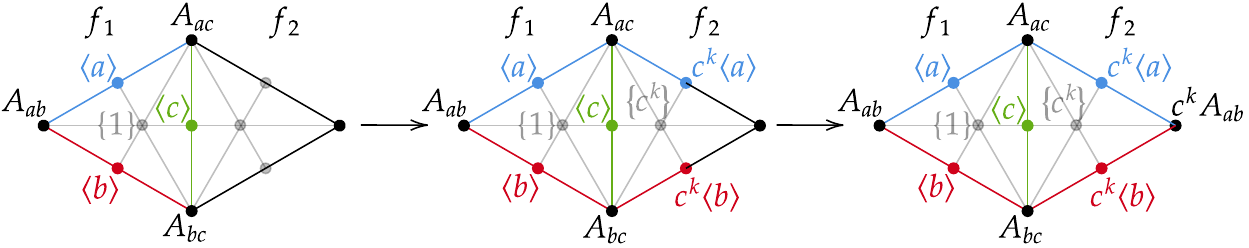}
\caption{Showing that $f_2$ is a translate of $f_1$. \underline{Left:} The initial situation. \underline{Middle:} Step 1. \underline{Right:} Step 2, from which we conclude Step 3.}
\label{FigTranslates}
\end{figure}

\noindent \underline{Proof of the claim:} We first assume that $f_1$ belongs to $D_2^i$. To prove the first assertion, we will want to think about the original thinner combinatorial structure on these triangles. Up to conjugation, we will suppose that $f_1$ is the triangle described in Figure \ref{FigTranslates}. We assign colours to the edges of $f_1$, giving two edges the same colour if and only if they have the same local group. We push these colours outside of $f_1$ by giving to any translate of an edge of $f_1$ the same colour as the edge it is the translate of. We want to show that the three edges of $f_2$ are of the same colours as those of $f_1$.
\\ \underline{Step 1:} $f_2$ is contained in a translate of the fundamental domain. Because $f_1$ and $f_2$ share the vertex $\langle c \rangle$, this translate has to be $c^k K_{\Gamma}$ for some $k \neq 0$. Thus the two type $1$ vertices of $f_2$ outside $f_1$ are $c^k \langle a \rangle$ and $c^k \langle b \rangle$.
\\ \underline{Step 2:} $f_2$ is a triangle, so there is a type $2$ vertex in $f_2$ that is a neighbour of both $c^k \langle a \rangle$ and $c^k \langle b \rangle$. The unique such type $2$ vertex is $c^k A_{ab}$. This allows to determine the colour (i.e. the orbit type) of every edge in $f_2$.
\\ \underline{Step 3:} The above directly implies that $f_2 = c^k f_1$.

\noindent Recall that $D_2^i$ is non-trivial. Hence using the above, one can show by induction starting at a triangle of $D_2^i$ that any two triangles of $D_2$ are translate of one another. This finishes the proof of the claim.
\medskip

\begin{figure}[H]
\centering
\includegraphics[scale=0.85]{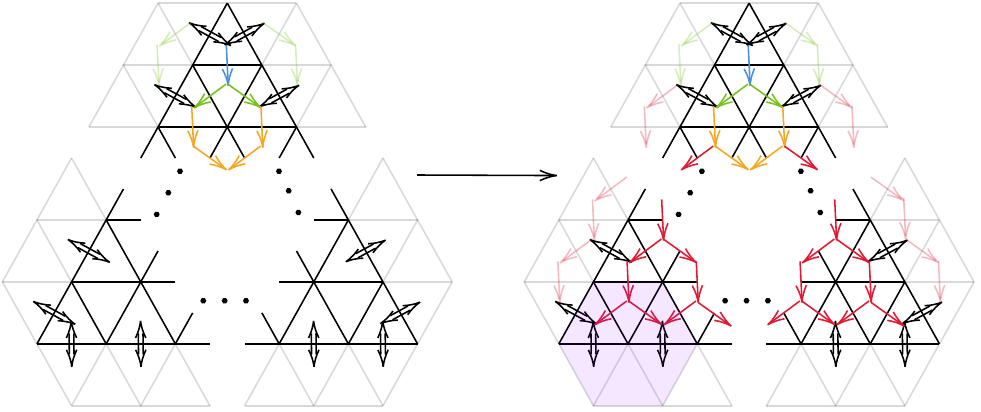}
\caption{Putting a system of arrows on $D$. \underline{Left:} Step 1 (black arrows), Step 2 (blue arrow), Step 3 (green arrows) and Step 4 (orange arrows). \underline{Right:} Step 5 (the induction process, red arrows). The purple hexagon gives a contradiction to $D$ being more than one triangle. The simplices and arrows not contained in $D$ are drawn with lighter colours.}
\label{FigM}
\end{figure}

Every triangle in $D_2^i$ is Euclidean equilateral. Hence using the claim, we actually obtain that every single triangle in $D$ is Euclidean equilateral. Thus the whole of $D$ is actually Euclidean, i.e. isometrically embedded in a Euclidean plane. A picture of $D$ is given in black on Figure \ref{FigM}.

Our next argument relies on the following notion:
\medskip

\noindent \textbf{Definition. (\cite[Definition 3.18]{vaskou2023isomorphism})} Let $g f$ and $h f$ be two adjacent equilateral triangles of $D$. Then there exists a standard generator $a \in V(\Gamma)$ and some integer $k \neq 0$ such that $g^{-1} h = a^k$. From that we define a \textbf{system of arrows} on the triangles of $D$:
\\(1) Draw a single arrow from $g f$ to $h f$ if $g^{-1} h = a$;
\\(2) Draw a double arrow between $g f$ and $h f$ if $g^{-1} h = a^k$ with $|k| \geq 2$.
\medskip

We now put a system of arrows on $D$. Consider a side $\gamma$ of $D$. We want to show that we can always put double arrows on the edges of $\gamma$. By hypothesis, and up to conjugation, $\gamma$ is contained in a standard tree $Fix(\langle a \rangle)$ for some $a \in V(\Gamma)$. We call the \textbf{strip} around $\gamma$ in $D$ the smallest subcomplex $S_{\gamma}$ of $D$ that contains every triangle of $D$ that has at least one vertex lying on $\gamma$. The strip $S_{\gamma}$ intersects its translate $a^2 S_{\gamma}$ exactly along $\gamma$, and the arrows between the simplices of $S_{\gamma}$ and that of $a^2 S_{\gamma}$ are all double arrows (see Figure \ref{FigM}).

As showed in (\cite[Lemma 3.20]{vaskou2023isomorphism}), systems of arrows are rather rigid, and must (up to symmetries or rotations) take one of the following forms:

\begin{figure}[H]
\centering
\includegraphics[scale=0.5]{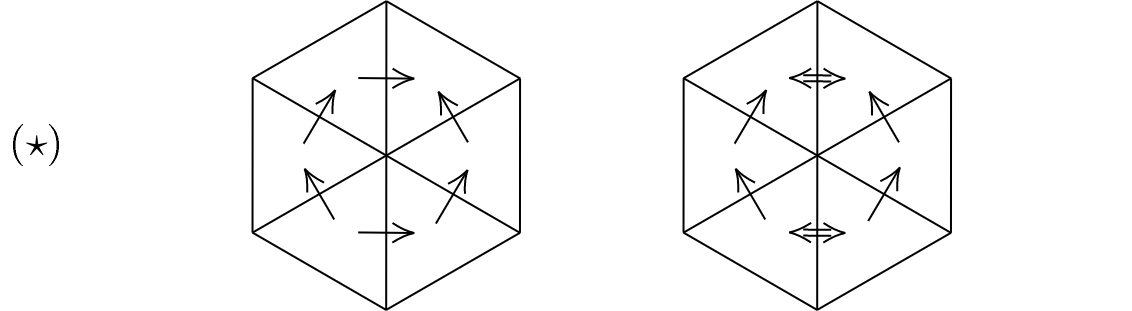}
\label{FigurePolygon333}
\end{figure}

\noindent In our case, because the arrows on the edges of $\gamma$ are all double arrows, the arrows between consecutive triangles of the strip $S_{\gamma}$ must all be single arrows. We now proceed to determine all the arrows in $D$:
\\ \underline{Step 1:} Put double arrows on the sides of $D$.
\\ \underline{Step 2:} The arrow between the two topmost triangles of $D$ must be simple by ($\star$). We suppose without loss of generality that it is pointing down.
\\ \underline{Step 3:} Use ($\star$) to complete the hexagons around this first arrow. We obtain two new arrows in $D$.
\\ \underline{Step 4:} Use ($\star$) on these two arrows and complete an hexagon of $D$.
\\ \underline{Step 5:} Proceed by induction using ($\star$) to determine every arrow in $D$.

We can see that the system of arrows of any of the hexagons along the bottommost side of $D$ contains two single arrows pointing away from each other and pointing towards double arrows (see Figure \ref{FigM}). This gives a contradiction to ($\star$). It follows that $D$ contains a single triangle. In particular, there is an element $g \in A_{\Gamma}$ such that $T \subseteq D \subseteq g K_{\Gamma}$, as wanted. Note that the element $g$ is unique with the property that $T \subseteq g K_{\Gamma}$, because $D$ is the unique embedded disc bounding $T$, and $g K_{\Gamma}$ is the unique translate of the fundamental domain that contains the vertex $\{g\}$, which belongs to $D$.
\hfill\(\Box\)
\bigskip

Before going to the next lemma we will need the following standard definition:

\begin{defi}
Let $X$ be a simplicial complex and let $x \in X$. Then the \textbf{star} of $x$ is the subcomplex $St(x)$ of $X$ spanned by $x$ and all the vertices adjacent to $x$. If $Y$ is a subcomplex of $X$, then the star of $x$ in $Y$ is the intersection $St_Y(x) \coloneqq St(x) \cap Y$.
\end{defi}

We recall that $\Gamma_{bar}$ can be seen as the boundary of the fundamental domain $K_{\Gamma}$, as explained in Remark \ref{RemGammaBar}.

\begin{lemma} \label{LemmaNeighOfType1} Let $x$ be a vertex of type $1$ of $X_{\Gamma}$, i.e. $x = g \langle a \rangle$ for some $g \in A_{\Gamma}$ and $a \in V(\Gamma)$. Then the star $St_{X_{\Gamma}^{(1)-ess}}(\langle a \rangle)$ of $x$ in $X_{\Gamma}^{(1)-ess}$ is the $g$-translate of the star $St_{\Gamma_{bar}}(\langle a \rangle)$ of $x$ in $\Gamma_{bar}$, and takes the form of a $n$-pod for some $n \geq 1$. It is contained in the standard tree $Fix(G_x)$, and in any translate of the fundamental domain that contains $x$.
\end{lemma}

\noindent \textbf{Proof:} First notice that $St_{X_{\Gamma}^{(1)-ess}}(x) = St_{X_{\Gamma}}(x) \cap X_{\Gamma}^{(1)-ess}$. By (\cite[II.12.24]{bridson2013metric}), the structure of $St_{X_{\Gamma}}(x)$ can be described as the development of a (sub)complex of groups that only depends on the local groups around $x$. Intersecting with $X_{\Gamma}^{(1)-ess}$ means further restricting to the local groups around $x$ that contain $G_x$. These local groups are the $g$-conjugates of the local groups around $\langle a \rangle$, so $St_{X_{\Gamma}^{(1)-ess}}(x)$ is the $g$-translate of $St_{\Gamma_{bar}}(\langle a \rangle)$, which is easily seen to be a $n$-pod, where $n$ is the number of edges attached to $\langle a \rangle$ in $\Gamma_{bar}$.

The inclusion $St_{X_{\Gamma}^{(1)-ess}}(x) \subseteq Fix(G_x)$ comes from the fact that every local group in the star contains $G_x$. Moreover, $St_{X_{\Gamma}^{(1)-ess}}(\langle a \rangle) \subseteq K_{\Gamma}$ and thus $St_{X_{\Gamma}^{(1)-ess}}(x) \subseteq h K_{\Gamma}$ for every $h \in A_{\Gamma}$ for which $x \in h K_{\Gamma}$.
\hfill\(\Box\)
\bigskip

\noindent \textbf{Proof of Proposition \ref{PropAdjEqAdj}:} [(1) $\Rightarrow$ (2)]: The vertices $v_1$ and $v_2$ are at combinatorial distance $2$ from each other, so there is a type $1$ vertex $x_{12}$ that is adjacent to both $v_1$ and $v_2$. Let us first suppose that $x_{12}$ belongs to $K_{\Gamma}$. By Lemma \ref{LemmaNeighOfType1}, $K_{\Gamma}$ contains the star $St_{X_{\Gamma}^{(1)-ess}}(x_{12})$, and this star is the simplicial neighbourhood of $x_{12}$ in $\Gamma_{bar}$. In particular then, $v_1$ and $v_2$ are distinct vertices of $\Gamma_{bar}$ that are adjacent to $x_{12}$. Because $\Gamma$ is complete, the path joining $v_1$, $x$ and $v_2$ can be completed into a cycle $\gamma \coloneqq (v_1, x_{12}, v_2, x_{23}, v_3, x_{31})$ of length $6$ in $\Gamma_{bar}$, where the $v_i$'s are type $2$ vertices and the $x_{ij}$'s are type $1$ vertices. Let now $H_3 \coloneqq f_{V_2}^{-1}(v_3)$. All that's left to do is to check that the pair $(H_1,H_2)$ satisfies the adjacency property, with respect to the third group $H_3$. This directly follows from Lemma \ref{LemmaEdgePropGeom}: the $H_i$'s are distinct subgroups, the sets $Fix(H_i \cap H_j)$ are distinct standard trees as they contain the type $1$ vertex $x_{ij}$ and no other type $1$ vertex of $\gamma$, and the trees $Fix(H_i \cap H_j)$ intersect pairwise along distinct type $2$ vertices, hence the triple intersection is trivial.
\medskip

\noindent If $x_{12}$ does not belong to $K_{\Gamma}$, then $x_{12} = g \cdot \bar{x}_{12}$, where $\bar{x}_{12}$ is a type $1$ vertex of $K_{\Gamma}$. Proceeding as before on $\bar{x}_{12}$ yields groups $H_i$ for $i \in \{1,2,3\}$. Then one can recover an analogous reasoning for $x_{12}$, using the groups $g H_i g^{-1}$ instead, for $i \in \{1,2,3\}$.

\medskip
\noindent [(2) $\Rightarrow$ (1)]: Let $(H_1, H_2)$ have the adjacency property relatively to a third subgroup $H_3$, and let $v_i \coloneqq f_{V_2}(H_i)$ for $i \in \{1, 2, 3\}$. By Lemma \ref{LemmaEdgePropGeom} the three $v_i$'s are distinct, and the three geodesics of the form $\gamma_{ij}$ connecting $v_i$ and $v_j$ are contained in the standard trees $Fix(H_i \cap H_j)$. The three $\gamma_{ij}$'s intersect pairwise, but the triple intersection is empty, by Lemma \ref{LemmaEdgePropGeom} again. Hence we are under the hypotheses of Lemma \ref{LemmaTrianglesInFundamentalDomains}, and we can conclude that the triangle formed by these three geodesics is contained is a single fundamental domain. In particular, the $v_i$'s are at combinatorial distance $2$ from each other.
\hfill\(\Box\)
\bigskip

\noindent We are now able to define explicitly the algebraic analogue of the type $1$ vertices of $X_{\Gamma}$:

\begin{defi} \label{DefiDV1} Let us consider the poset $\mathcal{P}_f(D_{V_2})$ of finite sets of distinct elements of $D_{V_2}$, ordered by the inclusion.
We now define $D_{V_1}$ to be the subset of $\mathcal{P}_f(D_{V_2})$ of sets $\{H_1, \cdots, H_k\}$ satisfying the following:
\medskip

\noindent \ \ \ \ \ \ \ \ \ \ \ \ \ \ \ (P1) \ Any subset $\{ H_i, H_j \} \subseteq \{H_1, \cdots, H_k\}$ is such that $(H_i, H_j)$

\noindent \ \ \ \ \ \ \ \ \ \ \ \ \ \ \ \ \ \ \ \ \ \ \ satisfies the adjacency property;
\smallskip

\noindent \ \ \ \ \ \ \ \ \ \ \ \ \ \ \ (P2) \ $\bigcap\limits_{i=1}^k H_i \neq \{1\}$;
\smallskip

\noindent \ \ \ \ \ \ \ \ \ \ \ \ \ \ \ (P3) \ $\{H_1, \cdots, H_k\}$ is maximal in $\mathcal{P}_f(D_{V_2})$ with these properties.
\end{defi}

\noindent As is was the case for the adjacency property, there is also a geometric meaning behind Definition \ref{DefiDV1}. While we managed to reconstruct the type $2$ vertices of $X_{\Gamma}$ directly from the spherical parabolic subgroups of type $2$ of $A_{\Gamma}$, we reconstruct a type $1$ vertex $x$ of $X_{\Gamma}$ from the sets of type $2$ vertices of $X_{\Gamma}$ that are adjacent to $x$. This is made more precise thereafter:

\begin{prop} \label{PropType1VertAlg1} The map $f_{V_1} : D_{V_1} \rightarrow V_1$ defined by the following is well-defined and is a bijection:
\\(1) For every element $\{H_1, \cdots, H_k \} \in D_{V_1}$, $f_{V_1}(\{H_1, \cdots, H_k \})$ is the unique vertex $x \in V_1$ that is adjacent to $v_i \coloneqq f_{V_2}(H_i)$ for every $H_i \in \{H_1, \cdots, H_k \}$.
\\(2) For every vertex $x \in V_1$, $f_{V_1}^{-1}(x)$ is the set $\{H_1, \cdots, H_k \} \in D_{V_1}$ of all the subgroups for which $v_i \coloneqq f_{V_2}(H_i)$ is adjacent to $x$.
\end{prop}

\noindent \textbf{Proof:} We first show that the two maps are well-defined. Then, checking that the composition of the two maps gives the identity is straightforward.
\medskip

\noindent \underline{$f_{V_1}$ is well-defined:} Let $\{H_1, \cdots, H_k \} \in D_{V_1}$. The intersection $H_1 \cap \cdots \cap H_k$ is an intersection of parabolic subgroups of type $2$ of $A_{\Gamma}$, hence is also a parabolic subgroup, by (\cite[Proposition 2.22.(1)]{vaskou2023isomorphism}). It is proper in any $H_i$ and non-trivial by definition, so it is a parabolic subgroup of type $1$ of $A_{\Gamma}$. The corresponding fixed set $T \coloneqq Fix(H_1 \cap \cdots \cap H_k)$ is a standard tree on which all the vertices $v_i \coloneqq f_{V_2}(H_i)$ lie. The convex hull $C$ of all the $v_i$'s in $T$ is a subtree of $T$. By hypothesis, any pair $(H_i, H_j)$ satisfies the adjacency property. Using Proposition \ref{PropAdjEqAdj}, this means the combinatorial distance between any two of the vertices defining the boundary of $C$ is $2$, so $C$ has combinatorial diameter $2$. As a tree with diameter $2$, $C$ contains exactly one vertex that is not a leaf of $C$, and this vertex must have type $1$.
\medskip

\noindent \underline{$f_{V_1}^{-1}$ is well-defined:} Let now $x \in V_1$, let $\{v_1, \cdots, v_k \}$ be the set of all the type $2$ vertices that are adjacent to $x$, and set $H_i \coloneqq f_{V_2}^{-1}(v_i)$. We want to check that $\{H_1, \cdots, H_k \} \in D_{V_1}$, i.e. that the properties (P1), (P2) and (P3) of Definition \ref{DefiDV1} are satisfied. First of all, we know that the combinatorial neighbourhood of $x$ is an $n$-pod that belongs to $Fix(G_x)$, by Lemma \ref{LemmaNeighOfType1}. In particular, all the $v_i$'s lie on the standard tree $Fix(G_x)$, which means that $G_x$ is contained in every $H_i$. This proves (P2).

Proving (P1) is straightforward if we use Proposition \ref{PropAdjEqAdj}: the $v_i$'s are distinct but they are all connected to a common vertex $x$, so the combinatorial distance between two distinct $v_i$'s is exactly $2$.

At last, if $\{H_1, \cdots, H_k \}$ was not maximal, there would be some $H_{k+1}$ such that $\{H_1, \cdots, H_{k+1}\}$ satisfies (P1) and (P2) of Definition \ref{DefiDV1}. The vertex $v_{k+1} \coloneqq f_{V_2}(H_{k+1})$ lies on $Fix(G_x)$ (use (P2)) and is at distance $2$ from all the other $v_i$'s (use (P1)), but is not adjacent to $x$ by hypothesis. This means one can connect $v_1$ and $v_2$ through $Fix(G_x)$ but without going through the star of $x$ in $Fix(G_x)$. This contradicts $Fix(G_x)$ being a tree. Therefore $\{H_1, \cdots, H_k \}$ is maximal, proving (P3).
\hfill\(\Box\)

\begin{rem} \label{RemAdjacent12} Let $H \in D_{V_2}$, let $\{H_1, \cdots, H_k \} \in D_{V_1}$, let $v \coloneqq f_{V_2}(H)$, and let $x \coloneqq f_{V_1}(\{H_1, \cdots, H_k \})$. Then one can easily deduce from the proof of Proposition \ref{PropType1VertAlg1} that $v$ and $x$ are adjacent if and only if $H \in \{H_1, \cdots, H_k \}$.
\end{rem}

\noindent We have reconstructed the algebraic analogue of the type $2$ vertices and the type $1$ vertices of $X_{\Gamma}$ (Lemma \ref{LemmaType2VertAlg} and Proposition \ref{PropType1VertAlg1}). To reconstruct the whole of $X_{\Gamma}^{(1)-ess}$, it remains to describe when an element of $D_{V_2}$ and an element of $D_{V_1}$ should be adjacent. This directly follows from Remark \ref{RemAdjacent12}:

\begin{defi} \label{DefiD1} We define the graph $D_1$ by the following:
\\(1) The vertex set of $D_1$ is the set $D_{V_2} \sqcup D_{V_1}$;
\\(2) We draw an edge between $H \in D_{V_2}$ and $\{H_1, \cdots, H_k \} \in D_{V_1}$ if and only if $H \in \{H_1, \cdots, H_k \}$.
\end{defi}

\begin{prop} \label{PropType1VertAlg} The bijections $f_{V_2}$ and $f_{V_1}$ can be extended to a graph isomorphism $F_1 : D_1 \rightarrow X_{\Gamma}^{(1)-ess}$.
\end{prop}

\noindent \textbf{Proof:} Let $f_{V_2} \sqcup f_{V_1} : D_{V_2} \sqcup D_{V_1} \rightarrow V_2 \sqcup V_1$. Then $f_{V_2} \sqcup f_{V_1}$ is a bijection by Lemma \ref{LemmaType2VertAlg} and Proposition \ref{PropType1VertAlg1}. We only need to show that two elements of $D_{V_2} \sqcup D_{V_1}$ are adjacent if and only if their images through $f_{V_2} \sqcup f_{V_1}$ are adjacent. Notice that
\begin{align*}
& H \in D_{V_2} \text{ and } \{H_1, \cdots, H_k \} \in D_{V_1} \text{ are adjacent in } D_1 \\
\overset{(\ref{DefiD1})} \Longleftrightarrow & H \in \{H_1, \cdots, H_k \} \\
\overset{(\ref{RemAdjacent12})} \Longleftrightarrow & f_{V_2}(H) \text{ and } f_{V_1}(\{H_1, \cdots, H_k \}) \text{ are adjacent in } X_{\Gamma}^{(1)-ess}.
\end{align*}
\hfill\(\Box\)

We have just reconstructed $X_{\Gamma}^{(1)-ess}$ purely algebraically. Our next goal is to reconstruct the whole of $X_{\Gamma}$. We start with the following definition:

\begin{defi} \label{DefiDV0} A subgraph $G$ of $D_1$ or of $X_{\Gamma}^{(1)-ess}$ is called \textbf{characteristic} if:
\medskip

\noindent \phantom{mmm} (C1) $G$ is isomorphic to the barycentric subdivision of a complete graph

\noindent \phantom{mmmmmm} on at least $3$ vertices;
\smallskip

\noindent \phantom{mmm} (C2) $G$ is maximal with that property.
\medskip

\noindent We call $\mathcal{CS}$ the set of characteristic subgraphs of $D_1$.
\end{defi}

\begin{lemma} \label{LemmaGSingleFD} The set of characteristic subgraphs of $X_{\Gamma}^{(1)-ess}$ is precisely $\{ g \Gamma_{bar} \ | \ g \in A_{\Gamma} \}$. In particular, $\mathcal{CS} = \{ F_1^{-1}(g \Gamma_{bar}) \ | \ g \in A_{\Gamma} \}$.
\end{lemma}

\noindent \textbf{Proof:} We focus on proving the first statement, as the second then directly follows from Proposition \ref{PropType1VertAlg}. We first prove the two following claims:
\medskip

\noindent \underline{Claim 1:} Any (non backtracking) cycle $\gamma \subseteq X_{\Gamma}^{(1)-ess}$ of length $6$ is contained in a single $g$-translate of the fundamental domain $K_{\Gamma}$.
\medskip

\noindent \underline{Proof of Claim 1:} Recall that $X_{\Gamma}^{(1)-ess}$ is a bipartite graph with partition sets $V_2$ and $V_1$. Consequently $\gamma = (x_1, v_{12}, x_2, v_{23}, x_3, v_{31})$, where the $x_i$'s are type $1$ vertices and the $v_{ij}$'s are type $2$ vertices of $X_{\Gamma}$. Consider now the three subgeodesics $c_1 \coloneqq (v_{31}, x_1, v_{12})$, $c_2 \coloneqq (v_{12}, x_2, v_{23})$ and $c_3 \coloneqq (v_{23}, x_3, v_{31})$, whose union is $\gamma$. Each geodesic $c_i$ is contained in the star $St_{X_{\Gamma}^{(1)-ess}}(x_i)$, which we know by Lemma \ref{LemmaNeighOfType1} is itself included in the standard tree $Fix(G_x)$. Also note that the three corresponding standard trees are distinct, or the fact that $\gamma$ is a cycle of length $6$ would contradict either the convexity of the standard trees, or the fact that they are uniquely geodesic. The three geodesics intersect pairwise, but their triple intersection is empty. We can now use Lemma \ref{LemmaTrianglesInFundamentalDomains}, and recover that $\gamma$ must be contained in a single translate $g K_{\Gamma}$. This finishes the proof of Claim 1.
\medskip

\noindent \underline{Claim 2:} For every subgraph $G$ of $X_{\Gamma}^{(1)-ess}$ that satisfies (C1) there exists an element $g \in A_{\Gamma}$ such that $G \subseteq g \Gamma_{bar}$.
\medskip

\noindent \underline{Proof of Claim 2:} $G$ is the barycentric subdivision of a complete graph $\widetilde{G}$ with at least $3$ vertices. In particular, $G$ contains at least one $6$-cycle, call it $\gamma_0$. This cycle corresponds to a $3$-cycle $\widetilde{\gamma_0}$ in $\widetilde{G}$. Because $\widetilde{G}$ is complete, every edge of $\widetilde{G}$ can be “reached” from $\widetilde{\gamma_0}$ by a sequence of $3$-cycles that consecutively intersect along an edge of $\widetilde{G}$. This means that for every edge $e$ of $G$, there exists a sequence of $6$-cycles $\gamma_0, \cdots, \gamma_n$ such that $e$ is contained in $\gamma_n$ and such that $\gamma_i, \gamma_{i+1}$ share exactly two edges (note that we can always pick $n \leq 2$). It is important to notice that the common vertex of these two edges has type $2$.

\begin{figure}[H]
\centering
\includegraphics[scale=0.85]{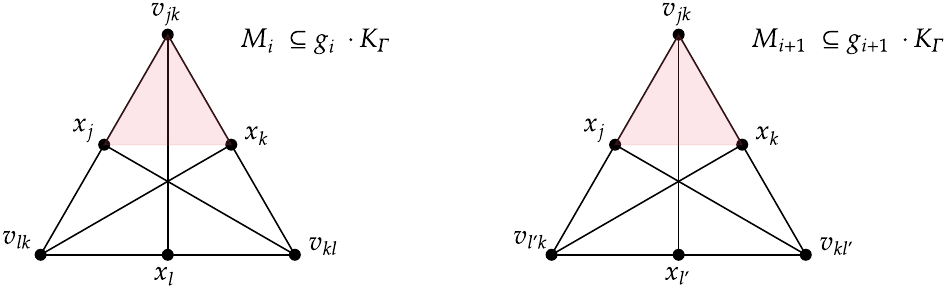}
\caption{The combinatorial subcomplexes $M_i$ (on the left) and $M_{i+1}$ (on the right). The convex hull $c(x_j, v_{jk}, x_k)$ is highlighted in light red, and is contained in both subcomplexes.}
\label{FigCommonFD}
\end{figure}

Using Claim 1, we know that each $\gamma_i$ is contained in a single translate $g_i K_{\Gamma}$. We want to show that all the $g_i$'s are the same element. To do so, we show that for every $0 \leq i < n$ we have $g_i = g_{i+1}$. Let $M_i$ be the subcomplex of $g_i K_{\Gamma}$ spanned by $\{g_i\}$ and by the vertices of $\gamma_i$. The two cycles $\gamma_i$ and $\gamma_{i+1}$ share two edges, whose union corresponds to a single edge of $\Gamma$. This means $M_i$ and $M_{i+1}$ share two edges of $X_{\Gamma}^{(1)-ess}$ that are attached to a common type $2$ vertex (see Figure \ref{FigCommonFD}). The convex hull of these two edges belongs to a single translate $g K_{\Gamma}$, yet belongs to both $g_i K_{\Gamma}$ and $g_{i+1} K_{\Gamma}$. This forces $g_i = g_{i+1}$. In particular, $e$ belongs to $g K_{\Gamma}$. As this works for every edge $e$ of $G$, we obtain $G \subseteq g K_{\Gamma}$.

\noindent Finally, $G$ is contained in the intersection $X_{\Gamma}^{(1)-ess} \cap g K_{\Gamma} = g \Gamma_{bar}$. This finishes the proof of Claim 2. We can now prove the main statement of the lemma:
\medskip

($\supseteq$) Consider a subgraph of $X_{\Gamma}^{(1)-ess}$ of the form $g \Gamma_{bar}$ for some $g \in A_{\Gamma}$. Because $\Gamma$ is a complete graph on at least $3$ vertices, $g \Gamma_{bar}$ checks (C1). If $g \Gamma_{bar}$ didn't satisfy (C2), there would be a graph $G$ that satisfies (C1) and strictly contains $g \Gamma_{bar}$. Using Claim 2, we obtain an element $g' \in A_{\Gamma}$ such that
$$g \Gamma_{bar} \subsetneq G \subseteq g' \Gamma_{bar}$$
By comparing the convex hulls as before, this forces $g = g'$, a contradiction.

\medskip
($\subseteq$) Let $G$ be a characteristic subgraph of $X_{\Gamma}^{(1)-ess}$. By Claim 2, there is an element $g \in A_{\Gamma}$ such that $G \subseteq g \Gamma_{bar}$. Using (C2) and the fact that $g \Gamma_{bar}$ is a characteristic subgraph shows this inclusion is an equality.
\hfill\(\Box\)

\begin{defi} \label{DefiAlgebraic} Let $D_{\Gamma}$ be the $2$-dimensional combinatorial complex defined by starting with $D_1$, and then coning-off every characteristic graph of $D_1$. The complex $D_{\Gamma}$ is called the \textbf{algebraic Deligne complex} associated with $A_{\Gamma}$.
\end{defi}

\begin{prop} \label{PropDIsomDeligne} The graph isomorphism $F_1$ from Proposition \ref{PropType1VertAlg} can be extended to a combinatorial isomorphism $F : D_{\Gamma} \rightarrow X_{\Gamma}$.
\end{prop}

\noindent \textbf{Proof:} We already know that the map $F_1$ of Proposition \ref{PropType1VertAlg} gives a graph isomorphism between $D_1$ and $X_{\Gamma}^{(1)-ess}$. The result now follows from the fact that $D_{\Gamma}$ and $X_{\Gamma}$ can respectively be obtained from $D_1$ and $X_{\Gamma}^{(1)-ess}$ by coning-off their characteristic subgraphs:
\\ $\bullet$ for $D_{\Gamma}$, this is the definition of the complex;
\\ $\bullet$ for $X_{\Gamma}$, this follows from Lemma \ref{LemmaGSingleFD} and Remark \ref{RemX1essCone}.
\hfill\(\Box\)
\bigskip

\noindent \textbf{Notation:} Following Proposition \ref{PropDIsomDeligne}, and to make the notation lighter, we will from now on slightly abuse the notation and identify $X_{\Gamma}$ with $D_{\Gamma}$, without caring about the combinatorial isomorphism $F$.
\bigskip

Theorem B as formulated in the introduction follows from the next theorem and its corollary, along with Proposition \ref{PropDIsomDeligne}.

\begin{thm} \label{ThmIsomorphicD} \textbf{(Theorem B)} Let $A_{\Gamma}$ and $A_{\Gamma'}$ be two large-type free-of-infinity Artin groups of rank at least $3$, with respective (algebraic) Deligne complexes $D_{\Gamma}$ and $D_{\Gamma'}$. Then any isomorphism $\varphi: A_{\Gamma} \rightarrow A_{\Gamma'}$ induces a natural combinatorial isomorphism $\varphi_* : D_{\Gamma} \rightarrow D_{\Gamma'}$, that can be described explicitly as follows:
\\ $\bullet$ For an element $H \in D_{V_2}^{\Gamma}$, $\varphi_*(H)$ is the subgroup $\varphi(H)$.
\\ $\bullet$ For a set $\{H_1, \cdots, H_k \} \in D_{V_1}^{\Gamma}$, $\varphi_*(\{H_1, \cdots, H_k \})$ is the set $\{ \varphi(H_1), \cdots, \varphi(H_k) \}$.
\\ $\bullet$ For an edge $e$ of $D_1^{\Gamma}$ connecting $H$ to $\{H_1, \cdots, H_k \}$, $\varphi_*(e)$ is the edge of $D_1^{\Gamma'}$ connecting $\varphi_*(H)$ to $\varphi_*(\{H_1, \cdots, H_k \})$.
\\ $\bullet$ For a simplex $f$ of $D_1^{\Gamma}$ connecting $H$, $\{H_1, \cdots, H_k \}$ and a vertex of type $0$ corresponding to the apex of a cone over a characteristic graph $G$, $\varphi_*(f)$ is the simplex of $D_1^{\Gamma'}$ connecting $\varphi_*(H)$, $\varphi_*(\{H_1, \cdots, H_k \})$, and the vertex of type $0$ corresponding to the apex of the cone over the characteristic graph $\varphi_*(G)$.
\end{thm}

\noindent \textbf{Proof:} The fact that $\varphi_*$ is a combinatorial isomorphism directly follows from the definition of $D_{\Gamma}$, that was constructed using algebraic notions such as inclusions, intersections and maximality, that are all preserved under isomorphisms. The explicit description of $\varphi_*$ is clear from the way we constructed the algebraic Deligne complexes.
\hfill\(\Box\)
\bigskip

A direct consequence of the previous theorem is the following corollary:

\begin{coro} \label{CoroCombiAction} There is a natural combinatorial action of $Aut(A_{\Gamma})$ onto $D_{\Gamma}$ (and $X_{\Gamma}$), and this action is explicitly described by the map $\varphi_*$ from Theorem \ref{ThmIsomorphicD}.
\end{coro}

\begin{rem} \label{RemAction}
The action of an automorphism $\varphi \in Aut(A_{\Gamma})$ on $D_{\Gamma}$ is entirely determined by its action on the set of type $2$ vertices of the complex. This is because every simplex of $D_{\Gamma}$, whether it is a type $1$ vertex, an edge, or a $2$-dimensional simplex, is defined algebraically from the set of type $2$ vertices of the complex.
\end{rem}

As a direct consequence of Theorem \ref{ThmIsomorphicD}, we show below that the class of large-type free-of-infinity Artin groups is rigid. Note that this result is already a consequence of (\cite[Theorem B]{vaskou2023isomorphism}).

\begin{thm} \label{TheoremRigidity}
Let $A_{\Gamma}$ and $A_{\Gamma'}$ be two large-type free-of-infinity Artin groups. Then $A_{\Gamma}$ and $A_{\Gamma'}$ are isomorphic as groups if and only if $\Gamma$ and $\Gamma'$ are isomorphic as labelled graphs.
\end{thm}

\noindent \textbf{Proof:} Consider an isomorphism $\varphi : A_{\Gamma} \rightarrow A_{\Gamma'}$. By Theorem \ref{ThmIsomorphicD}, $\varphi$ induces a combinatorial isomorphism $\varphi_* : D_{\Gamma} \rightarrow D_{\Gamma'}$ that sends the characteristic subgraphs of $D_1^{\Gamma}$ onto the characteristic subgraphs of $D_1^{\Gamma'}$. We claim that the isomorphism $\varphi_*: \Gamma_{bar} \rightarrow \varphi_*(\Gamma_{bar})$ is label-preserving. Indeed, every type $2$ vertex in $\Gamma_{bar}$ corresponds to a classical maximal dihedral Artin subgroup $H$ of $A_{\Gamma}$ with coefficient say $m$, and the corresponding type $2$ vertex in $\varphi_*(\Gamma_{bar})$ corresponds to the dihedral subgroup $\varphi(H)$, that also has coefficient $m$, because isomorphic dihedral Artin groups always have the same coefficients (see \cite[Theorem 1.1]{paris2003artin}). By Lemma \ref{LemmaGSingleFD}, there is an element $g \in A_{\Gamma'}$ such that $\varphi_*(\Gamma_{bar}) = g \Gamma'_{bar}$. Thus we have
$$\Gamma_{bar} \cong \varphi_*(\Gamma_{bar}) = g \Gamma_{bar}' \cong \Gamma_{bar}',$$
where every isomorphism is label-preserving. One then easily deduces that $\Gamma \cong \Gamma'$.
\hfill\(\Box\)

\section{Automorphism groups}

Let $A_{\Gamma}$ be a large-type free-of-infinity Artin group of rank at least $3$. In Section 3 we introduced various algebraic objects and proved that the Deligne complex $X_{\Gamma}$ associated with $A_{\Gamma}$ can be reconstructed in a purely algebraic way. This allowed to build a natural combinatorial action of the automorphism group $Aut(A_{\Gamma})$ onto the Deligne complex (Theorem B). In this section we will see that this action can be used to compute $Aut(A_{\Gamma})$ explicitly, proving Theorem A.

\begin{lemma} \label{LemmaInnIsAction} The group $Inn(A_{\Gamma})$ of inner automorphisms of $A_{\Gamma}$ acts on $X_{\Gamma}$ in a natural way: every inner automorphism $\varphi_g : h \mapsto g h g^{-1}$ acts on $X_{\Gamma}$ like the element $g$. Moreover $Inn(A_{\Gamma}) \cong A_{\Gamma}$.
\end{lemma}

\noindent \textbf{Proof:} We begin by proving the first statement. By Remark \ref{RemAction}.(1), it is enough to check that this holds when we restrict the action to type $2$ vertices of $X_{\Gamma}$. Let $g \in A_{\Gamma}$, and let $v \in V_2$ be a type $2$ vertex of $X_{\Gamma}$. Then
$$\varphi_g \cdot v \coloneqq (F \circ \varphi_g \circ F^{-1})(v) = F(\varphi_g (G_v)) = F(g G_v g^{-1}) = F(G_{g \cdot v}) = g \cdot v.$$
The fact that $Inn(A_{\Gamma}) \cong A_{\Gamma}$ is a consequence of $A_{\Gamma}$ having trivial centre (\cite[Corollary C]{vaskou2021acylindrical}).
\hfill\(\Box\)

\begin{defi} \cite[Definition 2.18]{vaskou2023isomorphism}
For any Artin group $A_{\Gamma}$, there is a well-defined homomorphism $ht: A_{\Gamma} \rightarrow \mathbf{Z}$ that sends every standard generator to $1$. For any element $g \in A_{\Gamma}$, we call $ht(g)$ the \textbf{height} of $g$.
\end{defi}

\begin{lemma} \label{LemmaHeightPreserving} Let $\iota$ be the automorphism of $A_{\Gamma}$ defined by $\iota(s) \coloneqq s^{-1}$ for every generator $s \in V(\Gamma)$, and let $\varphi \in Aut(A_{\Gamma})$ be any automorphism. Then one of $\varphi$ or $\varphi \circ \iota$ is height-preserving.
\end{lemma}

\noindent \textbf{Proof:} By Corollary \ref{CoroCombiAction} the automorphism $\varphi$ acts combinatorially on $X_{\Gamma}$. In particular, it sends the vertex $\{1\}$ onto the vertex $\{g\}$ for some $g \in A_{\Gamma}$. Using Lemma \ref{LemmaInnIsAction}, the automorphism $\varphi_{g^{-1}} \circ \varphi $ fixes $\{1\}$. Since inner automorphisms preserve height, we can simply suppose that $\varphi$ fixes $\{1\}$. In particular, $\varphi$ preserves $\Gamma_{bar}$ and thus sends the set of type $1$ vertices of $K_{\Gamma}$ onto itself. Looking at the action of $\varphi$ on $D_{\Gamma}$, this means $\varphi$ sends any standard parabolic subgroup of type $1$ of $A_{\Gamma}$ onto a similar subgroup. Consequently, every standard generator must be sent by $\varphi$ onto an element that generates such a subgroup, i.e. that has height $1$ or $-1$. There are three possibilities:
\medskip

\noindent \underline{(1) $ht(\varphi(s)) = 1, \forall s \in V(\Gamma)$:} Then $\varphi$ is height-preserving.
\smallskip

\noindent \underline{(2) $ht(\varphi(s)) = -1, \forall s \in V(\Gamma)$:} Then $\varphi \circ \iota$ is height-preserving.
\smallskip

\noindent \underline{(3) $\exists s, t \in V(\Gamma): ht(\varphi(s)) = 1$ and $ht(\varphi(t)) = -1$:} This means there are generators $a, b \in V(\Gamma)$ such that $\varphi(s) = a$ and $\varphi(t) = b^{-1}$. Because $A_{\Gamma}$ is free-of-infinity, the generators $s$ and $t$, as well as the generators $a$ and $b$, generate dihedral Artin subgroups of $A_{\Gamma}$. Note that $\varphi(A_{st}) = \langle \varphi(s), \varphi(t) \rangle = \langle a, b^{-1} \rangle = A_{ab}$. Because $\varphi$ is an isomorphism we must have $m_{st} = m_{ab}$ (use \cite[Theorem 1.1]{paris2003artin}). Applying $\varphi$ on both sides of the relation $sts \cdots = tst \cdots$ yields
$$ab^{-1}a \cdots = b^{-1}ab^{-1} \cdots.$$
Note that if we put everything on the same side, we obtain a word with $2 m_{st} = 2 m_{ab}$ syllables, that is trivial in $A_{ab}$. The words of length $2 m_{ab}$ that are trivial in $A_{ab}$ have been classified in (\cite[Lemma 3.1]{martin2020abelian}), and the word we obtained does not fit this classification, which yields a contradiction.
\hfill\(\Box\)

\begin{defi} Let $Aut(\Gamma)$ be the group of label-preserving graph automorphisms of $\Gamma$. We say that an isomorphism $\varphi \in Aut(A_{\Gamma})$ is \textbf{graph-induced} if there exists a graph automorphism $\phi \in Aut(\Gamma)$ such that $\varphi_*(\Gamma_{bar}) = \phi(\Gamma_{bar})$. We denote by $Aut_{GI}(A_{\Gamma})$ the subgroup of $Aut(A_{\Gamma})$ consisting of the graph-induced automorphisms.
\end{defi}

\begin{lemma} \label{LemmaGI1to1Aut} The map $\mathcal{F} : Aut_{GI}(A_{\Gamma}) \rightarrow Aut(\Gamma) \times \{id, \iota\}$ defined by the following is a group isomorphism:
\\Any $\varphi \in Aut_{GI}(A_{\Gamma})$ induces an automorphism of $\Gamma_{bar}$ and thus of $\Gamma$. This isomorphism defines the first component of $\mathcal{F}(\varphi)$. The second component of $\mathcal{F}(\varphi)$ is $id$ if $\varphi$ is height-preserving, and $\iota$ otherwise.
\end{lemma}

\noindent \textbf{Proof:} It is easy to check that $\mathcal{F}$ defines a morphism, so we show that it defines a bijection by describing its inverse map. Let $\phi \in Aut(\Gamma) \times \{id, \iota \}$. Then for any standard generator $s \in V(\Gamma)$, the automorphism $\phi$ sends the vertex $\langle s \rangle$ corresponding to $s$ onto the vertex $\phi(\langle s \rangle)$ corresponding to a standard generator that we note $s_{\phi}$. Define $\varphi_{\phi}$ as the (unique) automorphism of $A_{\Gamma}$ that sends every standard generator $s$ onto the standard generator $s_{\phi}$. Note that when acting on $X_{\Gamma}$, $\varphi_{\phi}$ restricts to an automorphism of $\Gamma_{bar}$ that corresponds to the automorphism $\varphi$ of $\Gamma$. For $\varepsilon \in \{0, 1 \}$ we let $\mathcal{F}^{-1}((\phi, \iota^{\varepsilon})) \coloneqq \varphi_{\phi} \circ \iota^{\varepsilon}$. It is clear that $\varphi_{\phi} \circ \iota^{\varepsilon}$ is graph-induced, and it is easy to check that composing $\mathcal{F}^{-1}$ with $\mathcal{F}$ on either side yields the identity.
\hfill\(\Box\)
\bigskip

We are now able to prove the main result of this paper:

\begin{thm} \label{MainThm} \textbf{(Theorem A)} Let $A_{\Gamma}$ be a large-type free-of-infinity Artin group of rank at least $3$. Then we have
$$Aut(A_{\Gamma}) \cong A_{\Gamma} \rtimes (Aut(\Gamma) \times \left(\quotient{\mathbf{Z}}{2 \mathbf{Z}})\right)
\ \ \ \text{ and } \ \ \ Out(A_{\Gamma}) \cong Aut(\Gamma) \times \left(\quotient{\mathbf{Z}}{2 \mathbf{Z}})\right).$$
\end{thm}

\noindent \textbf{Proof:} Let $\varphi \in Aut(A_{\Gamma})$. The same argument as the one in the proof of Lemma \ref{LemmaHeightPreserving} shows that up to post-composing with an inner automorphism, we may as well assume that $\varphi$ preserves $\Gamma_{bar}$, i.e. that $\varphi$ is graph-induced. This means
$$Aut(A_{\Gamma}) \cong Inn(A_{\Gamma}) \rtimes Aut_{GI}(A_{\Gamma}),$$
Using Lemma \ref{LemmaInnIsAction} and Lemma \ref{LemmaGI1to1Aut}, we obtain
$$Aut(A_{\Gamma}) \cong A_{\Gamma} \rtimes (Aut(\Gamma) \times \{id, \iota\}) \cong A_{\Gamma} \rtimes (Aut(\Gamma) \times \left(\quotient{\mathbf{Z}}{2 \mathbf{Z}})\right)).$$
In particular, we have
$$Out(A_{\Gamma}) \cong Aut(\Gamma) \times \left(\quotient{\mathbf{Z}}{2 \mathbf{Z}})\right).$$
\hfill\(\Box\)
\bigskip
\bigskip

\noindent \textbf{Acknowledgments:} I thank Alexandre Martin for our many discussions. This work was partially supported by the EPSRC New Investigator Award EP/S010963/1. I also thank the referee for their useful comments.

\nocite{*}

\newcommand{\etalchar}[1]{$^{#1}$}

\text{E-mail address:} \texttt{\href{mailto: nicolas.vaskou@bristol.ac.uk}{nicolas.vaskou@bristol.ac.uk}}

School of Mathematics, University of Bristol, Woodland Road Bristol BS8 1UG.

\end{document}